
\input amstex
\documentstyle{amsppt}
\magnification=1200
\catcode`\@=11
\redefine\logo@{}
\catcode`\@=13
\pageheight{19cm}

\define \bn{\Bbb N}

\define \bq{\Bbb Q}
\define \br{\Bbb R}
\define \bc{\Bbb C}

\define \M{{\Cal M}}
\define\Ha{{\Cal H}}
\define\La{{\Cal L}}
\define\E{{\Cal E}}

\define\rk{\text{rk}~}


\define\Exc{\text{Exc}}
\TagsOnRight
\NoBlackBoxes


\define\vertex#1{\hbox to15pt{\hfill \vbox to30pt{
\hbox to15pt{\hfill $#1$ \hfill}
\nointerlineskip\vfill
\hbox to15pt{\hfill $\bigcirc$ \hfill}
\nointerlineskip
\vskip10pt}\hfill}}

\define\vertexb#1{\hbox to20pt{\hfill \vbox to30pt{
\hbox to20pt{\hfill $#1$ \hfill}
\nointerlineskip\vfill
\hbox to20pt{\hfill $\bullet$ \hfill}
\nointerlineskip
\vskip13pt}\hfill}}

\define\darr{\hbox to30pt{\hfill
\vbox to30pt{\vfill \nointerlineskip
\hbox to25pt{\rightarrowfill} \nointerlineskip
\hbox to25pt{\leftarrowfill} \nointerlineskip
\vfill}
\hfill}}

\define\darrup#1{
\hbox to30pt{\hfill
\vbox to30pt{\vfill
\hbox to25pt{\hfill $#1$ \hfill}
\nointerlineskip
\hbox to25pt{\rightarrowfill} \nointerlineskip
\hbox to25pt{\leftarrowfill} \nointerlineskip
\vskip8pt}
\hfill}}

\define\darrdown#1{
\hbox to30pt{\hfill
\vbox to30pt{\vskip8pt \nointerlineskip
\hbox to25pt{\rightarrowfill} \nointerlineskip
\hbox to25pt{\leftarrowfill} \nointerlineskip
\hbox to25pt{\hfill $#1$ \hfill}
\vfill}
\hfill}}

\define\darrupdown#1#2{
\hbox to30pt{\hfill
\vbox to30pt{
\hbox to25pt{\hfill $#1$ \hfill} \nointerlineskip
\vskip2pt\vfill
\hbox to25pt{\rightarrowfill} \nointerlineskip
\hbox to25pt{\leftarrowfill} \nointerlineskip\vfill
\hbox to25pt{\hfill $#2$ \hfill}
\vfill}
\hfill}}

\define\horlines{\hbox to60pt{\hfil \vrule height30pt \hskip-3pt
\hbox to55pt{\hfil  \vbox to30pt{
\hrule width50pt
\nointerlineskip\vfil
\hrule width50pt
\nointerlineskip\vfil
\hrule width50pt
\nointerlineskip\vfil
\hrule width50pt
\nointerlineskip\vfil
\hrule width50pt
\nointerlineskip\vfil
\hrule width50pt
\nointerlineskip\vfil
\hrule width50pt
\nointerlineskip\vfil
\hrule width50pt
\nointerlineskip\vfil
\hrule width50pt
\nointerlineskip\vfil
\hrule width50pt
\nointerlineskip\vfil
\hrule width50pt}\hfil}
\hskip-6pt \vrule height30pt \hfil}}

\define\verlines{\hbox to60pt{\hfil
\vbox to30pt{
\hrule width50pt
\hbox to 50pt{\vrule height30pt
\hfil \vrule height30pt \hfil \vrule height30pt \hfil \vrule height30pt
\hfil \vrule height30pt \hfil \vrule height30pt \hfil \vrule height30pt
\hfil \vrule height30pt \hfil \vrule height30pt \hfil \vrule height30pt}
\hrule width50pt
}
\hfil}}

\document

\topmatter

\title
On algebraic varieties with finite polyhedral Mori cone 
\endtitle

\author
Viacheslav V. Nikulin \footnote{Supported by
Grant of Russian Fund of Fundamental Research (grant N 00-01-00170).
\hfill\hfill}
\endauthor

\address
Steklov Mathematical Institute,
ul. Gubkina 8, Moscow 117966, GSP-1, Russia
\endaddress
\address
Dept. of Pure Mathem. of the University of
Liverpool. Liverpool L69 3BX. England
\endaddress
\email
slava\@nikulin.mian.su \ \  V.Nikulin\@liverpool.ac.uk
\endemail

\abstract
The fundamental property of Fano varieties with mild singularities  
is that they have a finite polyhedral Mori cone. Thus, it is very 
interesting to ask: What we can say about algebraic varieties with a 
finite polyhedral Mori cone? I give a review of known results. 

All of them were obtained applying methods which were originated in 
the theory of discrete groups generated by reflections in hyperbolic 
spaces with a fundamental chamber of finite volume. 
\endabstract

\rightheadtext
{Algebraic varieties with finite polyhedral Mori cone}
\leftheadtext{V.V. Nikulin}

\endtopmatter

\head
1. Introduction
\endhead

Let $X$ be a projective algebraic variety with $\bq$-factorial 
singularities over an algebraically closed field. Let $N_1(X)$ be 
the $\br$-linear space generated by all algebraic curves on $X$ by 
the numerical equivalence, and the $N^1(X)$ be the $\br$-linear 
space generated by all Cartier (or Weil) divisors on $X$ by 
numerical equivalence. The $N_1(X)$ and $N^1(X)$ are dual to each other. 

The convex cone $NE(X)\subset N_1(X)$ generated by all effective curves 
on $X$ is called {\it Mori cone}. Its dual $NEF(X)=\{x\in N^1(X)\ |\ 
x\cdot NE(X)\ge 0$ is called {\it $nef$ cone.}    

One of basic properties of Fano varieties with log-terminal 
singularities which follows from Mori Theory \cite{Mo1}, \cite{Mo2} 
and its development by Kawamata \cite{Ka1}, and Shokurov \cite{Sho} 
is that the Mori cone $NE(X)$ is {\it finite polyhedral.}  
It is generated by a finite set of rays which are called 
{\it extremal rays}.  

Thus, it is interesting to ask: What can we say about algebraic 
varieties with finite polyhedral Mori cone? Another question could be: 
What can we get from this property for Fano varieties? 

Nowadays, we understand that well for surfaces, but there are 
some results also for 3-folds, for Fano 3-folds and Calabi--Yau 3-folds.  
All these results were obtained during last 15 years, and I want to make
a review of these results and also to speak about open problems. 

All methods, I know, were originated in the Theory of discrete groups 
generated by reflection in hyperbolic spaces with a fundamental 
chamber of finite volume. They were developed 
by the author and Vinberg more than 20 years ago. We exploit the idea 
that, in some cases, the polyhedron $\M(X)=(NEF(X)-\{0\})/\br^+$ 
has very similar properties to properties of a fundamental chamber 
$\M$ of finite volume of a discrete reflection group in a hyperbolic space. 
In some cases, Algebraic Geometry is very similar to Hyperbolic Geometry.

\head
2. Methods and results from the Theory of discrete reflection groups 
in hyperbolic spaces
\endhead

There are two methods which were developed in the Theory of discrete 
reflection groups $W$ generated by reflections in hyperbolic spaces $\La^n$. 
That is in simply-connected Riemannian manifolds $\La^n$ 
of a constant negative curvature, where $n\ge 2$ is dimension. They 
are based on studying of geometric and combinatorial properties of 
the fundamental chamber $\M$ of $W$. 

The first method is the {\it Method of narrow parts of polyhedra.} 
See \cite{N3}. It is a metric property of a finite closed 
convex polyhedron in a hyperbolic space. We formulate it in Sect. 
3 considering application of this method to algebraic surfaces.  
It is valid for any finite convex closed polyhedron in a hyperbolic space.  
Applying this method to fundamental chambers $\M$ of reflection 
groups, it was shown in \cite{N3} (see also \cite{N4}) that the number of 
maximal arithmetic reflection groups $W$ 
(equivalently, the number of polyhedra $\M$) is finite for the fixed 
dimension $n\ge 2$ of the hyperbolic space and the fixed degree 
$N=[K:\bq]$ of the ground field $K$ (the field of definition) of $W$. 
Here $K$ is some purely real algebraic number field. The  
number of fields $K$ of the fixed degree $N$ for all $n\ge 2$ 
is also finite. Arithmetic groups always have a finite fundamental 
domain of finite volume.   

Applying the Method of narrow parts of polyhedra to algebraic varieties, 
we expect to get boundedness, existence of a very ample divisor of bounded 
degree. We consider these results in Sect. 3. The method of 
narrow parts of polyhedra is now developed and applied only 
to surfaces. It would be very interesting to develop it for 3-folds.  

The second method is the {\it Diagram Method.} It was started for 
arithmetic reflection groups in \cite{N4} and was developed for arbitrary   
reflection groups with bounded fundamental chambers $\M$ by 
Vinberg in \cite{V}. It uses  the combinatorial 
property of a bounded fundamental chamber $\M$ which is that 
the number of faces of highest dimension of $\M$ containing a 
vertex is equal to dimension. Equivalently, the dual polyhedron 
to $\M$ is simplicial. The important development of this method was 
obtained by Khovansky \cite{Kh} and Prokhorov \cite{P}. They generalized 
the method for polyhedra $\M$ which are simplicial in 1-dimensional 
faces. This is important for reflection groups with fundamental chambers 
$\M$ of finite volume but with some vertices at infinity. We formulate 
the method in Sect. 4 considering algebraic surfaces. Using this method, 
it was proved in \cite{N4} for arithmetic reflection groups 
that the dimension $n\le 9$ if the degree 
$N$ is large; in \cite{V} that $n\le 29$ for arbitrary 
discrete reflection groups with bounded fundamental chamber; in 
\cite{Kh}, \cite{P} that  $n\le 995$ for arbitrary discrete reflection 
groups with fundamental chamber of finite volume.               

Applying the Diagram Method to algebraic varieties, we expect to 
get an estimate for the Picard number $\rho=\dim N_1(X)$ which is absolute 
for the considered class of algebraic varieties. 
This method is very effective for Del Pezzo surfaces with log-terminal  
singularities. See Sect. 4. It is very interesting that there are some  
generalizations of this method also to 3-fold: Fano 3-folds and 
Calabi--Yau 3-folds. See Sect. 5.

\head
3. The method of narrow parts of polyhedra and algebraic surfaces.  
\endhead

We consider non-singular projective algebraic surfaces $X$ over 
an algebraically closed field, and we assume that $\rho=\dim N_1(X)\ge 3$. 
If the Mori cone $NE(X)$ is finite polyhedral, 
by Hodge Index Theorem, 
$$
NE(X)=\sum_{E\in \Exc(X)}{{\br}^+E}
$$
is generated by the finite set $\Exc(X)$ of all exceptional curves $E$ on 
$X$. We remind that a curve $E\subset X$ is {\it exceptional,} if $E$ is 
irreducible and $E^2<0$. We have natural invariants 
$$
\rho=\dim N_1(X)\ge  3,\ 
\delta=\max_{E\in \Exc(X)} {(-E^2)},\
p=\max_{E\in \Exc(X)} {(p_a(E))} 
$$
where $p_a(E)={E^2+E\cdot K\over 2}+1$ is the arithmetic genus of $E$ and 
$K=K_X$ the canonical class.  

We have 

\proclaim{Theorem 1 (\cite{N15})} 
For any algebraic surface $X$ with finite polyhedral Mori cone   
and the invariants $\rho\ge 3$, $\delta$, $p$,  
there exists a very ample divisor $H$ such that 
$$
H^2\le N(\rho,\delta,p)
$$
where $N(\rho,\delta,p)$ is a constant depending only on 
$\rho$, $\delta$ and  $p$. 
\endproclaim

There are examples in \cite{N15} that Theorem 1 is not true if at least 
one of constants $\rho$, $\delta$, $p$ is not fixed. 
The proof of Theorem 1 
follows from standard results on symmetric matrices with 
non-negative coefficients (the Perron--Frobenius Theorem) and 
the following result. 

\proclaim{Lemma 1 (Method of narrow parts of polyhedra 
\cite{N3, Appendix})}
There exist  
$E_1,...,E_\rho\in Exc(X)$ such that 

1) $E_1,...,E_\rho$ generate $N_1(X)$;

2) ${4(E_i\cdot E_j)^2\over E_i^2\,E_j^2}<62^2$  
(an absolute constant) for any $1\le i,j\le n$ ;

3) The Dynkin diagram of $E_1,...,E_\rho$ is connected (the set cannot be 
divided into two non-empty orthogonal to each other subsets). 
\endproclaim 

Using Lemma 1, the very ample divisor $H$ of Theorem 1 can be obtained as 
$H=a_1E_1+\cdots +a_\rho E_\rho$ where $a_i\in \bn$ depend 
only on the intersection matrix 
$(E_i\cdot E_j)$, $1\le i,j\le \rho$,  
the number of these matrices is finite.   

\demo{Idea of Proof of Lemma 1} 
It is very geometrical. I can explain how one can find the first element 
$E_1$. 

Let us fix  a {\it very ample} $D$. We find $E_1, F\in Exc(X)$ such
$[E_1, F]$ is negative definite,
$(F\cdot D)/\sqrt{-F^2}+(E_1\cdot D)/\sqrt{-E_1^2}$ is maximal,  and
$(F\cdot D)/\sqrt{-F^2}\ge
(E_1\cdot D)/\sqrt{-E_1^2}$. Roughly speaking,
we are looking for an exceptional curve $F$ of maximal degree with 
respect to $D$.  

By the Hodge Index Theorem, the space $N_1(X)$ is hyperbolic for the 
intersection pairing. It has the 
signature $(1,\rho-1)$. Then one can relate with $N_1(X)$ the 
hyperbolic space $\La(X)=V^+(X)/\br^+$ where 
$V^+(X)$ is the half (containing $D$) of the {\it light cone} $V(X)$ of all 
elements $x\in N_1(X)$ with $x^2>0$. The $\La (X)$ is equipped with 
the hyperbolic distance  
$|\br^+x\br^+y|=(x,y)/\sqrt{x^2y^2}$.  
The polyhedron 
$$
\M(X)=NEF(X)/\br^{+}=\{\br^+x\in \La(X)\ |\ x\cdot Exc(X)\ge 0\}\subset \La(X)
$$
is a finite closed convex polyhedron in $\La (X)$. The set  
$Exc(X)\subset N_1(X)$ gives all orthogonal vectors to faces of $\M(X)$ of 
codimension one. An element $E\in Exc(X)$ gives the   
face of codimension one of $\M(X)$ containing in the hyperplanes  
$\Ha_ E=\{\br^+x\ |\ x\cdot E=0\}\subset \La (X)$ which is orthogonal to $E$.

The same Lemma 1 is 
valid for any finite closed convex polyhedron $\M$ in $\La (X)$ if 
one replaces $Exc(X)$ by the set $P(\M)\subset N_1(X)$ of orthogonal 
vectors to faces of $\M$ of codimension one. Lemma 1 does not change if 
one replaces $E_i$ by $\lambda_iE_i$, $\lambda_i>0$. The statement 2)  
means that the hyperplanes of the corresponding faces are on absolutely 
bounded distance to each other. The statements 1) and 3) mean that 
they are in {\it general position.}  
\enddemo

Because of Theorem 1, it is interesting to classify surfaces $X$ with 
finite polyhedral Mori cone for small invariants 
$\rho$, $\delta$ and $p$. It is done for 
$\delta=2$ and $p=0$, see \cite{N15} for details. 
Equivalently, $E^2=-1$ or $-2$, and $E$ is 
non-singular rational for any exceptional curve $E\in \Exc(X)$.  
Then any $E\in Exc(X)$ defines a reflection $x\to x-2(x\cdot E)/E^2$, 
$x\in NS(X)$. Here $NS(X)\subset N_1(X)$ is the Neron-Severi lattice of 
$X$ generated by all algebraic curves on $X$. Thus, $\M(X)$ is 
a fundamental chamber of the corresponding arithmetic reflection group, 
and then $\rho=\dim N_1(X)$ is absolutely bounded. The surfaces $X$ are 
divided in several types (a) --- (d) described below.  

\smallpagebreak 

(a) {\it Minimal resolutions of singularities of Del Pezzo surfaces with 
Du Val (or index 1 log-terminal) singularities.} 
This is the case $K^2>0$ for the canonical 
class $K$ of the surface. All these surfaces are rational, $\rho \le 9$. 
Classification of Del Pezzo surfaces with Du Val 
singularities was obtained by Nagata in \cite{Na}. It is also a part of 
classification of Del Pezzo surfaces with log-terminal singularities of 
index $\le 2$ which was obtained in \cite{AN1} and \cite{AN2}. We shall 
discuss it in Sect. 4.

\smallpagebreak     

(b) {\it Rational surfaces with $K^2=0$ and finite polyhedral Mori cone.} 
This is the case $K^2=0$, but $K\not\equiv 0$. Then $\rho=10$. 
Classification of these surfaces was obtained in \cite{N15}. There are two 
different cases of that surfaces. First case: $\dim |-nK|=0$ for any 
$n\ge 1$. These surfaces $X$ are classified by the graphs of exceptional 
curves $\Exc(X)$, there are three possible graphs, given below (the black 
vertex denotes an exceptional curve $E$ with $E^2=-2$).  

\smallpagebreak 

\centerline{\hbox {
\vbox{\hbox{$E_{10}$} \hbox{$\circ$}}
\hskip-10pt
\hbox to1cm{\hrulefill}
\hskip-1pt
\vbox{\hbox{$E_9$} \hbox{$\bullet$}}
\hskip-6pt
\hbox to1cm{\hrulefill}
\hskip-1pt
\vbox{\hbox{$E_8$} \hbox{$\bullet$}}
\hskip-6pt
\hbox to1cm{\hrulefill}
\hskip-1pt
\vbox{\hbox{$E_7$} \hbox{$\bullet$}}
\hskip-6pt
\hbox to1cm{\hrulefill}
\hskip-1pt
\vbox{\hbox{$E_6$} \hbox{$\bullet$}}
\hskip-6pt
\hbox to1cm{\hrulefill}
\hskip-1pt
\vbox{\hbox{$E_5$} \hbox{$\bullet$}}
\hskip-6pt
\hbox to1cm{\hrulefill}
\hskip-1pt
\hbox{\vbox to 0.8cm{\vskip0.15cm\hbox{$E_4$} \hbox{$\bullet$}
\vskip1pt
\hbox{\hskip2pt \vrule height 0.8cm width 0.5pt} \vskip-5pt\hbox{$\bullet$}
\hbox{$E_1$}}}
\hskip-6pt
\hbox to1cm{\hrulefill}
\hskip-1pt
\vbox{\hbox{$E_3$} \hbox{$\bullet$}}
\hskip-6pt
\hbox to1cm{\hrulefill}
\hskip-1pt
\vbox{\hbox{$E_2$} \hbox{$\bullet$}}
}}

\vskip2cm

\centerline{Graph $H\widetilde{E_8}$}

\vskip1cm

\centerline{\hbox {
\vbox{\hbox{$E_{11}$} \hbox{$\circ$}}
\hskip-10pt
\hbox to1cm{\hrulefill}
\hskip-1pt
\vbox{\hbox{$E_{10}$} \hbox{$\bullet$}}
\hskip-6pt
\hbox to1cm{\hrulefill}
\hskip-1pt
\hbox{\vbox to 0.8cm{\vskip0.15cm\hbox{$E_7$} \hbox{$\bullet$}
\vskip1pt
\hbox{\hskip2pt \vrule height 0.8cm width 0.5pt} \vskip-5pt\hbox{$\bullet$}
\hbox{$E_1$}}}
\hskip-6pt
\hbox to1cm{\hrulefill}
\hskip-1pt
\vbox{\hbox{$E_5$} \hbox{$\bullet$}}
\hskip-6pt
\hbox to1cm{\hrulefill}
\hskip-1pt
\vbox{\hbox{$E_3$} \hbox{$\bullet$}}
\hskip-6pt
\hbox to1cm{\hrulefill}
\hskip-1pt
\vbox{\hbox{$E_2$} \hbox{$\bullet$}}
\hskip-6pt
\hbox to1cm{\hrulefill}
\hskip-1pt
\hbox{\vbox to 0.8cm{\vskip0.15cm\hbox{$E_6$} \hbox{$\bullet$}
\vskip1pt
\hbox{\hskip2pt \vrule height 0.8cm width 0.5pt} \vskip-5pt\hbox{$\bullet$}
\hbox{$E_4$}}}
\hskip-6pt
\hbox to1cm{\hrulefill}
\hskip-1pt
\vbox{\hbox{$E_8$} \hbox{$\bullet$}}
\hskip-6pt
\hbox to1cm{\hrulefill}
\hskip-1pt
\vbox{\hbox{$E_9$} \hbox{$\circ$}}
}}

\vskip1.5cm

\centerline{Graph $H\widetilde{D_8}$}

\vskip1cm

\centerline{
\hbox to8cm {\hfill
\vbox{\hbox{$\circ$} \hbox{$E_{12}$}}
\hskip-10pt
\vbox to 1.9cm{\hbox{\hskip0.6cm$E_1$}\vfill
\hbox to1cm{\hrulefill}\vskip15pt}
\hskip-5pt
\hbox{
\vbox{\hbox{$E_4$}
\hbox{$\bullet$}
\hbox{\hskip2pt \vrule height 0.8cm width 0.5pt}
\hbox{$\bullet$}
\hbox{\hskip2pt \vrule height 0.8cm width 0.5pt}
\hbox{$\bullet$}
\hbox{$E_9$}}
\hskip-6pt
\vbox{\hbox to1cm{\hrulefill}\vskip15pt}
\hskip-1pt
\vbox to 3cm{
\hbox{\hskip-1.1cm \hbox to1.6cm{\hrulefill}}
\vfill
\hbox{$\bullet$} \hbox{$E_6$}}
\hskip-9pt
\vbox to 4.5cm{
\hbox{\hskip11pt $E_{10}$}
\hbox{\hskip11pt $\circ$}
\hbox{\hskip11pt \hskip2pt \vrule height 0.6cm width 0.5pt}
\vskip-4pt
\hbox{\hskip11pt $\bullet$}
\hbox{\hskip11pt $E_7$}
\vfill
\hbox to1.1cm{\hrulefill}\vskip15pt}
\hskip-1pt
\vbox to 3cm{
\hbox{\hskip-0.4cm \hbox to1.8cm{\hrulefill}}
\vfill
\hbox{$\bullet$} \hbox{$E_3$}}
\hskip-33pt
\vbox{\hbox to1cm{\hrulefill}\vskip15pt}
\hskip-3pt
\hbox{
\vbox{
\hbox{$E_2$}
\hbox{$\bullet$}
\hbox{\hskip2pt \vrule height 0.8cm width 0.5pt}
\hbox{$\bullet$}
\hbox{\hskip2pt \vrule height 0.8cm width 0.5pt}
\hbox{$\bullet$}
\hbox{$E_8$}
}}
\hskip-11pt
\hbox{
\vbox to 1.9cm{
\hbox{$E_5$}
\vfill
\hbox to1cm{\hrulefill}\vskip15pt
}
\hskip-2pt
\vbox{
\hbox{$\circ$} \hbox{$E_{11}$}
}
}}\hfill}
}

\vskip1cm

\centerline{Graph $H\widetilde{A_8}$}

\smallpagebreak 

Second case: $\dim |-nK|=1$ for some $n\ge 1$. The linear 
system $|-nK|$ defines then an elliptic fibration on $X$ with a singular fibre 
of multiplicity $n$. 
The invariant $n\in \bn$ can take any value, and 
the surfaces $X$ have infinite number of connected components of moduli. 
Construction  and classification of these surfaces uses general 
theory of rational elliptic surfaces \cite{H}, \cite{CD}, \cite{D1},  
and Ogg-Shafarevich Theory of elliptic surfaces \cite{O}, \cite{Sha}. 
To have a finite polyhedral Mori cone, the sum of ranks of degenerate fibres 
of the elliptic fibration $|-nK|$ must be maximal, equals to $8$. 
Equivalently, the automorphism group of $X$ 
(it is mainly the Mordell--Weil group) must be finite.  

\smallpagebreak 

(c) {\it K3 surfaces with $\rho \ge 3$ and finite automorphism group.}  
This is the case when $K=0$. All exceptional curves $E$ have $E^2=-2$. 
By the description of automorphism groups of K3 surfaces 
in \cite{P-\u S\u S}, for $\rho\ge 3$
the Mori cone $NE(X)$ is finite polyhedral, if and only if 
$Aut (X)$ is finite (over $\bc$). All K3 surfaces $X$ with $\rho \ge 3$ and 
finite automorphism group were classified in 
\cite{N1}, \cite{N2}, \cite{N5} and \cite{N7}. For $\rho=4$, 
this classification is due to Vinberg, see the result in \cite{N7}.   
These K3 surfaces are classified by their Picard lattices $S=NS(X)$ 
which have the property that the group $W^{(-2)}(S)$ generated by 
reflections in all roots $\alpha$ of $S$ with $\alpha^2=-2$ has finite 
index in the automorphism group $O(S)$ of the lattice. There is a finite 
number of such hyperbolic lattices $S$ of $\rho=\rk (S)\ge 3$.  
All of them were found in \cite{N1}, \cite{N2}, \cite{N5} and \cite{N7}.
 Here is their number for each $\rho \ge 3$: 
$$
\matrix
3& 4 & 5 & 6 & 7 & 8 & 9 & 10 & 11 & 12& 13 & 
14&15&16&17&18&19& \ge 20 \\ 
27&17&10 & 10& 9 &12 & 10& 9 & 4 &  4  & 3 &  
3 & 1 & 1 & 1&  1& 1& 0
\endmatrix 
$$
A K3 surface $X$ may have very arbitrary Picard lattice 
$S$. For example, any even hyperbolic lattice $S$ of rank $\le 11$ is 
a Picard lattice of some K3 surface. But only finite number of $S$ with 
$\rho =\rk S\ge 3$ corresponds to K3 surfaces with finite 
automorphism group, and equivalently with finite polyhedral Mori cone. 
K3 surfaces with $\rho \ge 3$ having a finite polyhedral Mori cone 
are extremely rare and exceptional.

\smallpagebreak

(d) {\it Enriques surfaces with finite automorphism group.} This is 
the case when $K\equiv 0$, but $K\not=0$. All exceptional curves 
$E$ have $E^2=-2$, the $\rho=10$.  There are 
seven types of Enriques surfaces with finite automorphism group. Two of them 
depend on 1 moduli, and five give isolated Enriques surfaces. They
were classified in \cite{N6} (see also \cite{N7}), and by Kond\B o 
in \cite{Ko} where one missing in \cite{N6} type was also found.  

The first example of Enriques surfaces with finite automorphism group was 
found by Gino Fano \cite{F} in 1910. It is hard to follow Fano, but his 
example really gives an Enriques surface with the finite automorphism group! 
Since Enriques surfaces depend on 10 parameters, Fano's ideas    
were genius and correct. This paper by Fano was rediscovered by 
Dolgachev in 1985 during his visit of the University of Torino where 
he looked Fano's archive. One year before, 
Dolgachev published the paper \cite{D2} where he discovered  
the first contemporary example of Enriques surfaces depending on 
one parameter with finite automorphism group. It is completely non-obvious 
that Enriques surfaces with finite automorphism group do exist! 
The papers \cite{N6}, \cite{N7} and \cite{Ko} followed this paper.   

\smallpagebreak 

Unfortunately, no 3-dimensional generalization of the Method of 
narrow parts of polyhedra is known, in spite of its very clear  
geometric ideas.

\head
4.The Diagram method for surfaces 
\endhead

Let $X$ be a projective algebraic surface with isolated singularities and 
$\pi:\widetilde{X}\to X$ the minimal resolution of singularities.
By results of  Mori \cite{Mo1} applied to non-singular surfaces, 
$\widetilde{X}$ has a {\it finite polyhedral Mori cone,} if 
$\rho(\widetilde{X}) \ge 3$. The Mori cone $NE(\widetilde{X})$ is  
generated by exceptional curves which are in pre-images  
of singular points, and by exceptional curves of the first kind on 
$\widetilde{X}$. 

The diagram method can be applied to all these surfaces $\widetilde{X}$,  
and all surfaces with finite polyhedral Mori cone as well, but  
we get the most astonishing results if singularities are log-terminal (or 
$X$ is a log Del Pezzo surface).  
Over $\bc$ log-terminal singularities are the same as 
quotient singularities. 

Applying Diagram Method, we get 

\proclaim{Theorem 2 (\cite{N8}---\cite{N10})} There are functions $A(n)$, 
$B(n)$ of $n\in \bn$ such that for 
any Del Pezzo surface $X$ with log-terminal singularities one has:  
$$
\rho(\widetilde{X})<A(maximal\ index\ of\ singularities);
$$
$$
\rho(\widetilde{X})<B(maximal\ multiplicity\ of\ singularities).
$$
It follows boundedness for moduli of $X$ if maximal index or multiplicity 
of singularities is bounded (e.g. apply Theorem 1). 
\endproclaim

\proclaim{Theorem 3 (Alexeev, \cite{A})} There is an absolute 
constant $C>0$ such that for any Del Pezzo surface $X$ with log-terminal 
singularities 
$$
\rho(\widetilde{X})<{C\over \epsilon(X)}\ \ \ \text{where}
$$
$$
1\ge \epsilon(X)=\min_{E\,of\,1st\,kind} {(-E\cdot \pi^\ast K_X)}>0.
$$
It follows (Shokurov's Conjecture) that the set of fractional indices of 
all log Del Pezzo surfaces is contained in $(0,1]$, 
all its accumulation points are 
$1/n$, $n\in \bn$, and all of them are accumulation points 
from above, but not from below.
\endproclaim 

Remind that a positive rational number $r>0$ is {\it the fractional index} 
of a log Del Pezzo surface $X$ if $K_X=rH$ where $H\in Pic X$ is 
a primitive element of the Picard lattice.

Proofs of Theorems 2, 3 are based on the Lemma 2 below which constitutes 
the Diagram Method. A subset ${\Cal L}\subset \Exc(X)$ is called 
{\it Lanner} if it is minimal hyperbolic. A subset  ${\Cal E}\subset \Exc(X)$ 
is {\it elliptic,} if it is negative definite. The set $\Exc (X)$ 
defines a graph in a usual way. The set of ifs vertices is 
$\Exc (X)$. Two vertices $E_1$, $E_2$ are connect by an edge if 
$E_1\cdot E_2>0$. Thus, we can consider the distance, the diameter and so on.

\proclaim{Lemma 2 (Diagram Method, \cite{N8}, \cite{N14})} 
Let $\widetilde{X}$ be an algebraic surface with a finite 
polyhderal Mori cone. Assume that
 
(a) $\text{diam\ }(\La)  \le d$ for any Lanner subset $\La 
\subset \Exc (\widetilde{X})$;

(b)
$$
\#\{ \{ E_1,E_2\} \subset \E \ \mid \ 1 \le \rho (E_1,E_2)\le d \} \le
C_1\# \E
$$
and
$$
\#\{ \{ E_1,E_2\} \subset \E \ \mid \ d+1\le \rho (E_1,E_2)\le 2d+1\} \le
C_2\# \E
$$
for any elliptic subset $\E \subset Exc (\widetilde{X})$. 
Then
$$
\rho(\widetilde{X}) < 96(C_1+C_2/3)+68.
$$
\endproclaim

To prove Theorems 2 and 3, one should just show that the constants 
$d$, $C_1$ and $C_2$ of Lemma 2 can be estimated by    
functions depending on the maximal index, maximal multiplicity of 
singularities of $X$, or the invariant $1/\epsilon (X)$. 
This is not easy, but it can be done. 
See \cite{N8}---\cite{N11} and \cite{A}.

\demo{Idea of Proof of Lemma 2} Like for the proof of Lemma 1, we have 
a finite closed hyperbolic polyhedron    
$$
\M=NEF(\widetilde{X})/\br^{+}=
\{\br^+x\in \La(\widetilde{X})\ |\ x\cdot Exc(\widetilde{X})\ge 0\}
\subset \La({\widetilde{X}})
$$
of dimension $n=\rho-1$ with the set $\Exc(\widetilde{X})\subset N_1(X)$ 
of orthogonal vectors to faces of 
$\M$ of codimension one. Since $E\cdot E^\prime\ge 0$ for 
$E\not=E^\prime\in  \Exc (\widetilde{X})$, 
the polyhedron $\M$ has acute angles. By Perron-Frobenius Theorem, 
$\M$ is simplicial in its 
finite vertices and in its edges (1-dimensional faces).  

For simplicity, let us assume that all vertices of $\M$ are finite. 
Then $\M$ is dual to a simplicial polyhedron. It follows  that the 
combinatorial polynomial of $\M$ 
(where $\alpha_i$ is the number of $i$-dimensional faces of $\M$),   
$$
R(s)=\alpha_0+(s-1)\alpha_1+\cdots +\alpha_{n-1}(s-1)^{n-1}+(s-1)^n, 
$$
is reversible and has positive coefficients (e.g. see \cite{St}). 
It follows \cite{N4} that the average number 
$A^{0,2}$ of vertices of plane faces of $\M$ is bounded as 
$A^{0,2}\le  4+4/(n-2)$. Thus, almost all 2-dimensional faces of 
$\M$ are quadrangles or triangles if $n$ is large. In particular, since 
$\M$ is a hyperbolic polyhedron, $\M$ has a lot of non-right plane angles.   

A plane angle $A$ of $\M$ is defined by an elliptic subset 
$\E\subset \Exc(\widetilde{X})$ with $n$ vertices and by two its 
distinguished elements $E_1\not=E_2\in \E$. The vertex of $A$ is defined as 
intersection of all hyperplane orthogonal to elements of $\E$. Two  
sides of the angle $A$ are intersections of all hyperplanes which are 
orthogonal to elements of $\E-\{E_1\}$ and $\E-\{E_2\}$. 
The angle $A$ is called {\it combinatorially right, } 
if the distance between $E_1$ and $E_2$ in  
$\E$ is greater than $2d+1$. Otherwise, $A$ is not combinatorially right.

It is easy to see that a Lanner subset $\La$ is connected. Two Lanner subsets 
cannot be orthogonal. It follows that any triangle of $\M$ has at least two 
combinatorially non-right angles. Any quadrangle of $\M$ has at least one 
combinatorially non-right angle. Since almost all 2-dimensional faces of $\M$ 
are triangles or quadrangles, $\M$ has a lot of combinatorially non-right 
angles. On the other hand, by the condition (b) of of the lemma, 
$\M$ has not too many of combinatorially non-right angles. 
If $n$ is big, we get a contradiction, which gives some estimate on 
dimension of $\M$ depending on the constants $d$, $C_1$ and $C_2$. 

In \cite{N4}, for arithmetic reflection groups, usual right 
angles in hyperbolic geometry were used. 
Introducing of the combinatorial notion of 
right angles is due to Vinberg \cite{V}. The generalization of these 
considerations for polyhedra $\M$ with some vertices at infinity is due 
to Khovanskii \cite{Kh} and Prokhorov \cite{P}. They had to use 
3-dimensional angles of $\M$ instead of 2-dimensional. Already in  
\cite{N4},  3-dimensional angles were used with the estimate 
$A^{2,3}\le 6+12/(n-2)$ for the average number of 2-dimensional 
faces in $3$-dimensional faces of $\M$. In \cite{N4}, a general estimate for 
the average number or $k$-dimensional faces in $m$-dimensional faces of 
$\M$ was given. Khovanskii \cite{Kh}, without using of combinatorial 
polynomials, had shown that the same estimates are valid 
for polyhedra $\M$ which are simplicial only in edges.     

See details in \cite{N14} and the cited papers. 
\enddemo 

Theorem 2 shows that for a fixed $k$, it is possible in principle to classify 
Del Pezzo surfaces with log-terminal singularities of index $\le k$. 
It is now known only for $k=2$. 
In \cite{AN1} and \cite{AN2}, classification of Del Pezzo surfaces $X$  
with log-terminal singularities of index $\le 2$ is given. It is shown that 
the linear system $|-2K_X|$ has a non-singular curve $C$ which 
does not contain singular points of $X$. There exists an 
appropriate right resolution of singularities $p:X_1\to X$ such that 
$|-2K_{X_1}|$ has a non-singular curve 
$\widetilde{C}=p^{-1}(C)+E_1+\cdots +E_k$ where $E_1,\dots , E_k$ 
are non-singular rational curves with $E_i^2=-4$. 
The double covering of $X_1$ ramified in $\widetilde{C}$ gives a K3 
surface $Y$ with a non-symplectic involution $\sigma$ of the double 
covering. It reduces classification of $X$ to classification 
of K3 surfaces $Y$ with non-symplectic involutions $\sigma$ 
such that the fixed point set $Y^\sigma$ has a connected 
component of genus $\ge 2$. The surface $X_1$ has a finite polyhedral 
Mori cone with the polyhedron  $\M(X_1)=NEF(X_1)/\br^+$ which is 
a fundamental chamber for an arithmetic reflection group. 
The classification of the Del Pezzo surfaces 
$X$ enumerates possible graphs of the finite sets 
$\Exc(X_1)$ of exceptional curves on the surfaces $X_1$. 
There are plenty of cases.     

\smallpagebreak 
  
In \cite{K-M} a different geometric approach to log Del Pezzo surfaces 
is developed. It is based on studying of rational curves on the log 
Del Pezzo surfaces.  

Boundedness of log Del Pezzo surfaces of a fixed index (Theorem 2 for 
index) is now generalized by Alexandr Borisov \cite{B} 
to Fano 3-folds with log-terminal singularities. Similar result   
for $n$-dimensional Fano varieties was recently announced by 
McKernan \cite{McK}. All other results of Theorems 2 and 3 are now  
known only for surfaces, and their proofs use the diagram method.   

We can say that {\it the property to have a finite 
polyhedral Mori cone for $\widetilde{X}$ is now crucial for description  
of Del Pezzo surfaces $X$ with log-terminal singularities.}

\head
5. The Diagram Method for 3-folds
\endhead 

Generalizing the diagram method, two results about 3-folds were 
obtained. One of them is valid for Fano 3-folds with terminal 
$\bq$-factorial singularities, another one is valid for  
3-dimensional Calabi--Yau manifolds. 

We remind that an extremal ray $R$ of the Mori cone $NE(X)$ is called 
{\it divisorial} if effective curves of $R$ fill out a divisor $D(R)$ 
of $X$; in theorems below it is always irreducible. An extremal ray  
$R$ is called {\it small} if effective curves of $R$ fill out a curve of $X$. 
A proper face $\gamma$ of $NE(X)$ has {\it Kodaira dimension 
$d$, $d\le 3$,} if the contraction of $\gamma$ 
(in theorems below, it aways exists) is a morphism such that its image 
has dimension $d$. An irreducible curve $C$ belongs to $\gamma$ if and only if 
its image is a point for the contraction of $\gamma$.

\proclaim{Theorem 4 (\cite{N12})} 
Let $X$ be a Fano 3-fold with terminal
{\bf Q}-factorial singularities (by \cite{Mo1}, \cite{Mo2}, it has 
a finite polyhedral Mori cone).

Then
$$
\rho (X)\le 7
$$
except the following two cases:

(1) The Mori cone $NE(X)$ has a face of Kodaira dimension
$\le 2$ (i.e. its contraction gives a fibration
$\pi:X\to Y$ where $\dim Y\le 2$).

(2) There exists a small extremal ray (it does not exist if 
$X$ is non-singular).
\endproclaim

Of course, we know a lot about non-singular Fano 3-folds and 
about Fano 3-folds with terminal singularities, e. g. see \cite{M-M} and 
\cite{Ka2}, and this result is not something extraordinary. But 
Theorem 4 is obtained by a very elementary, purely combinatorial method. 

\smallpagebreak 

We have the following statement about 3-dimensional Calabi--Yau 
manifolds. 

\proclaim{Theorem 5 (Nikulin and Shokurov, \cite{N13})} 
Let $X$ be a 3-dimensional Calabi--Yau manifold 
(here the Mori cone is not necessarily finite polyhedral, and we have to add 
the case (3) below).

Then
$$
\rho (X)\le 40
$$
except the following three cases:

(1) There exists a rational $nef$ element $D$ with $D^3=0$ 
(it is expected that $|nD|$ gives a fibration of $X$ for big $n$). 

(2) There exists a small extremal ray 
(the corresponding flop gives then another birational model).

(3) The Mori cone is not finite polyhedral (it is expected that 
$Aut (X)$ is then infinite, \cite{Mor}). 
\endproclaim

\demo{Sketch of Proofs} If we exclude exceptions (1) --- (3),  
Mori cone of $X$ is finite polyhedral,  
all extremal rays are divisorial, there are no $nef$ rational 
elements $D$ with $D^3=0$. Further proof is divided in several steps. 

{\bf Step A:} We prove that two different extremal rays $R_1$ and $R_2$ 
cannot have the same divisor $D(R_1)=D(R_2)$. The picture below

\vskip0.5cm

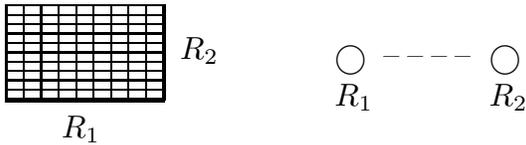
\line{\hskip10pt
\vbox{\horlines \hbox to60pt{\hfil $R_1$ \hfil}}
\hskip-62pt
\vbox{\verlines\hbox to60pt{\hfil}}
\hskip-5pt
\raise25pt\hbox{$R_2$}
\hskip30pt
\raise10pt\hbox{
\vbox{\vskip-10pt \hbox{$\bigcirc$} \hbox{$R_1$}}
\hskip0pt
\raise13pt\hbox to30pt{\hfil -- -- -- -- \hfil }
\hskip0pt
\vbox{\vskip-10pt \hbox{$\bigcirc$} \hbox{$R_2$}}}
\hfil}

\vskip0.2cm

\noindent 
is impossible. This is one of the most difficult steps which reflects problems
of 3-dimensional geometry. The proof for Calabi--Yau manifolds  
is especially complicated.

\vskip0.3cm

{\bf Step B:} We can {\it introduce diagrams} like below:

\vskip10pt

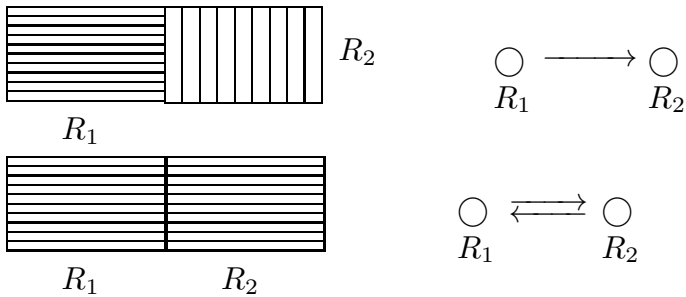
\line{\hskip10pt
\vbox{\horlines \hbox to60pt{\hfil $R_1$ \hfil}}
\hskip-12pt
\vbox{\verlines\hbox to60pt{\hfil}}
\hskip-5pt
\raise25pt\hbox{$R_2$}
\hskip30pt
\raise10pt\hbox{
\vbox{\vskip-10pt \hbox{$\bigcirc$} \hbox{$R_1$}}
\hskip0pt
\raise13pt\hbox to30pt{\rightarrowfill}
\hskip0pt
\vbox{\vskip-10pt \hbox{$\bigcirc$} \hbox{$R_2$}}}
\hfil}

\vskip0.1cm

\line{\hskip10pt
\vbox{\horlines \hbox to60pt{\hfil $R_1$ \hfil}}
\hskip-13pt
\vbox{\horlines\hbox to60pt{\hfil $R_2$ \hfil}}
\hskip30pt
\raise10pt\hbox{
\vbox{\vskip-10pt \hbox{$\bigcirc$} \hbox{$R_1$}}
\hskip0pt
\darr
\hskip0pt
\vbox{\vskip-10pt \hbox{$\bigcirc$} \hbox{$R_2$}}}
\hfil}

\vskip0.2cm 

Here arrow $R_1\to R_2$ means $R_1\cdot D(R_2)>0$. Thus, 
the set $\Exc (X)$ of all extremal rays (all of them are 
divisorial) defines an oriented graph. Using this graph, we 
can consider distances in subsets of $\Exc (X)$, their diameters.  
Since the graph is oriented, the distance function might be not symmetric. 

\vskip0.2cm

{\bf Step C:} 
A subset $\E\subset \Exc (X)$ is called 
{\it elliptic} if it is contained in a proper face of $NE(X)$. A subset
$\La\subset \Exc(X)$ is called a {\it E-set (extremal)} 
if it is minimal non-elliptic. (E-sets are similar to Lanner subsets 
for surfaces.) We have the following statement which  
constitutes the diagram method for 3-folds.

\proclaim{Lemma 3 (Diagram Method, \cite{N12}, \cite{N13})}
Assume
that there are 
constants $d$, $C_1$, $C_2$ such that the conditions (a) and (b)
below hold:

(a)
$$
\text{diam}(\La )\le d
$$
for any $E$-subset $\La\subset \Exc (X)$. 

(b)
$$
\sharp \{ (R_1, R_2)\in \E \times \E
\mid 1 \le \rho (R_1,R_2)\le d\} \le C_1 \sharp \E ;
$$
and
$$
\sharp \{ (R_1, R_2)\in \E \times \E
\mid d+1\le \rho (R_1,R_2) \le 2d+1\} \le C_2 \sharp \E .
$$
for any extremal subset $\E\subset \Exc(X)$. 

Then 
$$
\rho(X)  \le (16/3)C_1+4C_2+6.
$$
\endproclaim

\demo{About proof of Lemma 3} 
The proof of Lemma 3 is very similar to the proof of Lemma 2 for surfaces. 
Again one should work with the polyhedron $\M(X)=NEF(X)/\br^+$. 
Its faces are also orthogonal to elements of $\Exc(X)$. 
The crucial step is to prove that $\M(X)$ is dual to a  
simplicial polyhedron. Then one can use the same estimates for  
the average number of vertices of 2-dimensional faces of $\M(X)$. 
Further, considering {\it oriented plane angles} of $\M(X)$ instead of 
non-oriented, one should go through all steps of the proof of Lemma 2 
with some small changes.   
\enddemo   

{\bf Step D:} To apply Lemma 3, one needs to describe (classify) elliptic 
subsets and E-subsets. For 3-dimensional Calabi-Yau manifolds, one can 
see some of them below:

\vskip0.3cm

\centerline{\hfil
\hbox to 245pt{
\hfil
\hbox to15pt{\hfill \vbox to30pt{\vfill ${\bold A}_n\ :$ \vfill}\hfill}
\hskip10pt plus2pt minus2pt
\hbox to15pt{\hfill \vbox to30pt{\vfill $\bigcirc$ \vfill}\hfill}
\hskip8pt plus2pt minus2pt
\hbox to30pt{\hfill
\vbox to30pt{\vfill
\nointerlineskip
\hbox to25pt{\rightarrowfill} \nointerlineskip
\hbox to25pt{\leftarrowfill} \nointerlineskip
\vfill}
\hfill}
\hskip-14pt plus2pt minus2pt
\hbox to15pt{\hfill \vbox to30pt{\vfill $\bigcirc$ \vfill}\hfill}
\hskip8pt plus2pt minus2pt
\hbox to30pt{\hfill
\vbox to30pt{\vfill \nointerlineskip
\hbox to25pt{\rightarrowfill} \nointerlineskip
\hbox to25pt{\leftarrowfill} \nointerlineskip
\vfill}
\hfill}
\hskip-14pt plus2pt minus2pt
\hbox to15pt{\hfill \vbox to30pt{\vfill $\bigcirc$ \vfill}\hfill}
\hskip5pt plus2pt minus2pt
\hbox to15pt{\hfill \vbox to30pt{\vfill $\cdots$ \vfill}\hfill}
\hskip5pt plus2pt minus2pt
\hbox to15pt{\hfill \vbox to30pt{\vfill $\bigcirc$ \vfill}\hfill}
\hskip8pt plus2pt minus2pt
\hbox to30pt{\hfill
\vbox to30pt{\vfill \nointerlineskip
\hbox to25pt{\rightarrowfill} \nointerlineskip
\hbox to25pt{\leftarrowfill} \nointerlineskip
\vfill}
\hfill}
\hskip-14pt plus2pt minus2pt
\hbox to15pt{\hfill \vbox to30pt{\vfill $\bigcirc$ \vfill}\hfill}
\hskip8pt plus2pt minus2pt
\hbox to30pt{\hfill
\vbox to30pt{\vfill \nointerlineskip
\hbox to25pt{\rightarrowfill} \nointerlineskip
\hbox to25pt{\leftarrowfill} \nointerlineskip
\vfill}
\hfill}
\hskip-14pt plus2pt minus2pt
\hbox to15pt{\hfill \vbox to30pt{\vfill $\bigcirc$ \vfill}\hfill}
\hfil
}
\hskip55pt}


\centerline{\hfil
\hbox to 245pt{
\hfil
\hbox to15pt{\hfill \vbox to30pt{\vfill ${\bold B}_n\ :$ \vfill}\hfill}
\hskip10pt plus2pt minus2pt
\hbox to15pt{\hfill \vbox to30pt{\vfill $\bigcirc$ \vfill}\hfill}
\hskip8pt plus2pt minus2pt
\hbox to30pt{\hfill
\vbox to30pt{\vskip8pt \nointerlineskip
\hbox to25pt{\rightarrowfill} \nointerlineskip
\hbox to25pt{\leftarrowfill} \nointerlineskip
\hbox to25pt{\hfill $2$ \hfill}
\vfill}
\hfill}
\hskip-14pt plus2pt minus2pt
\hbox to15pt{\hfill \vbox to30pt{\vfill $\bigcirc$ \vfill}\hfill}
\hskip8pt plus2pt minus2pt
\hbox to30pt{\hfill
\vbox to30pt{\vfill \nointerlineskip
\hbox to25pt{\rightarrowfill} \nointerlineskip
\hbox to25pt{\leftarrowfill} \nointerlineskip
\vfill}
\hfill}
\hskip-14pt plus2pt minus2pt
\hbox to15pt{\hfill \vbox to30pt{\vfill $\bigcirc$ \vfill}\hfill}
\hskip5pt plus2pt minus2pt
\hbox to15pt{\hfill \vbox to30pt{\vfill $\cdots$ \vfill}\hfill}
\hskip5pt plus2pt minus2pt
\hbox to15pt{\hfill \vbox to30pt{\vfill $\bigcirc$ \vfill}\hfill}
\hskip8pt plus2pt minus2pt
\hbox to30pt{\hfill
\vbox to30pt{\vfill \nointerlineskip
\hbox to25pt{\rightarrowfill} \nointerlineskip
\hbox to25pt{\leftarrowfill} \nointerlineskip
\vfill}
\hfill}
\hskip-14pt plus2pt minus2pt
\hbox to15pt{\hfill \vbox to30pt{\vfill $\bigcirc$ \vfill}\hfill}
\hskip8pt plus2pt minus2pt
\hbox to30pt{\hfill
\vbox to30pt{\vfill \nointerlineskip
\hbox to25pt{\rightarrowfill} \nointerlineskip
\hbox to25pt{\leftarrowfill} \nointerlineskip
\vfill}
\hfill}
\hskip-14pt plus2pt minus2pt
\hbox to15pt{\hfill \vbox to30pt{\vfill $\bigcirc$ \vfill}\hfill}
\hfil
}
\hskip55pt}


\centerline{\hfil
\hbox to 245pt{
\hfil
\hbox to15pt{\hfill \vbox to30pt{\vfill ${\bold C}_n\ :$ \vfill}\hfill}
\hskip10pt plus2pt minus2pt
\hbox to15pt{\hfill \vbox to30pt{\vfill $\bigcirc$ \vfill}\hfill}
\hskip8pt plus2pt minus2pt
\hbox to30pt{\hfill
\vbox to30pt{\vfill
\hbox to25pt{\hfill $2$ \hfill}
\nointerlineskip
\hbox to25pt{\rightarrowfill} \nointerlineskip
\hbox to25pt{\leftarrowfill} \nointerlineskip
\vskip8pt}
\hfill}
\hskip-14pt plus2pt minus2pt
\hbox to15pt{\hfill \vbox to30pt{\vfill $\bigcirc$ \vfill}\hfill}
\hskip8pt plus2pt minus2pt
\hbox to30pt{\hfill
\vbox to30pt{\vfill \nointerlineskip
\hbox to25pt{\rightarrowfill} \nointerlineskip
\hbox to25pt{\leftarrowfill} \nointerlineskip
\vfill}
\hfill}
\hskip-14pt plus2pt minus2pt
\hbox to15pt{\hfill \vbox to30pt{\vfill $\bigcirc$ \vfill}\hfill}
\hskip5pt plus2pt minus2pt
\hbox to15pt{\hfill \vbox to30pt{\vfill $\cdots$ \vfill}\hfill}
\hskip5pt plus2pt minus2pt
\hbox to15pt{\hfill \vbox to30pt{\vfill $\bigcirc$ \vfill}\hfill}
\hskip8pt plus2pt minus2pt
\hbox to30pt{\hfill
\vbox to30pt{\vfill \nointerlineskip
\hbox to25pt{\rightarrowfill} \nointerlineskip
\hbox to25pt{\leftarrowfill} \nointerlineskip
\vfill}
\hfill}
\hskip-14pt plus2pt minus2pt
\hbox to15pt{\hfill \vbox to30pt{\vfill $\bigcirc$ \vfill}\hfill}
\hskip8pt plus2pt minus2pt
\hbox to30pt{\hfill
\vbox to30pt{\vfill \nointerlineskip
\hbox to25pt{\rightarrowfill} \nointerlineskip
\hbox to25pt{\leftarrowfill} \nointerlineskip
\vfill}
\hfill}
\hskip-14pt plus2pt minus2pt
\hbox to15pt{\hfill \vbox to30pt{\vfill $\bigcirc$ \vfill}\hfill}
\hfil
}
\hskip55pt}


\centerline{\hfil
\hbox to 290pt{
\hfil
\hbox to15pt{\hfill \vbox to60pt{\vfill ${\bold D}_n\ :$ \vfill}\hfill}
\hskip25pt plus2pt minus2pt
\hbox to 20pt{\hfill
\vbox to 60pt {\vfill
\hbox to 15pt{\hfill $\bigcirc$ \hfill}
\nointerlineskip\vfill
\hbox to15pt{\hfill $\uparrow\downarrow$ \hfill}
\nointerlineskip\vfill
\hbox to15pt{\hfill $\bigcirc$ \hfill}
\nointerlineskip\vfill
\hbox to15pt{\hfill $\uparrow\downarrow$ \hfill}
\nointerlineskip\vfill
\hbox to15pt{\hfill $\bigcirc$ \hfill}
\vfill}
\hfill}
\hskip-10pt plus2pt minus2pt
\hbox to30pt{\hfill
\vbox to60pt{\vfill
\nointerlineskip
\hbox to25pt{\rightarrowfill} \nointerlineskip
\hbox to25pt{\leftarrowfill} \nointerlineskip
\vfill}
\hfill}
\hskip-14pt plus2pt minus2pt
\hbox to15pt{\hfill \vbox to60pt{\vfill $\bigcirc$ \vfill}\hfill}
\hskip5pt plus2pt minus2pt
\hbox to15pt{\hfill \vbox to60pt{\vfill $\cdots$ \vfill}\hfill}
\hskip5pt plus2pt minus2pt
\hbox to15pt{\hfill \vbox to60pt{\vfill $\bigcirc$ \vfill}\hfill}
\hskip8pt plus2pt minus2pt
\hbox to30pt{\hfill
\vbox to60pt{\vfill \nointerlineskip
\hbox to25pt{\rightarrowfill} \nointerlineskip
\hbox to25pt{\leftarrowfill} \nointerlineskip
\vfill}
\hfill}
\hskip-14pt plus2pt minus2pt
\hbox to15pt{\hfill \vbox to60pt{\vfill $\bigcirc$ \vfill}\hfill}
\hskip8pt plus2pt minus2pt
\hbox to30pt{\hfill
\vbox to60pt{\vfill \nointerlineskip
\hbox to25pt{\rightarrowfill} \nointerlineskip
\hbox to25pt{\leftarrowfill} \nointerlineskip
\vfill}
\hfill}
\hskip-14pt plus2pt minus2pt
\hbox to15pt{\hfill \vbox to60pt{\vfill $\bigcirc$ \vfill}\hfill}
\hfil
}
\hskip50pt}


\centerline{\hfil
\hbox to 220pt{
\hfil
\hbox to15pt{\hfill \vbox to30pt{\vskip2pt${\bold E}_6\ :$ \vfill}\hfill}
\hskip10pt plus2pt minus2pt
\hbox to15pt{\hfill\vbox to30pt{\vskip2pt $\bigcirc$ \vfill}\hfill}
\hskip15pt plus2pt minus2pt
\hbox to30pt{\hfill
\vbox to30pt{
\nointerlineskip
\hbox to25pt{\rightarrowfill} \nointerlineskip
\hbox to25pt{\leftarrowfill} \nointerlineskip
\vfill}
\hfill}
\hskip-14pt plus2pt minus2pt
\hbox to15pt{\hfill \vbox to30pt{\vskip2pt $\bigcirc$ \vfill}\hfill}
\hskip8pt plus2pt minus2pt
\hbox to30pt{\hfill
\vbox to30pt{\nointerlineskip
\hbox to25pt{\rightarrowfill} \nointerlineskip
\hbox to25pt{\leftarrowfill} \nointerlineskip
\vfill}
\hfill}
\hskip-5pt plus2pt minus2pt
\hbox to 20pt{\hfill
\vbox to 30pt {\vskip2pt
\hbox to15pt{\hfill $\bigcirc$ \hfill}
\nointerlineskip\vfill
\hbox to15pt{\hfill $\uparrow\downarrow$ \hfill}
\nointerlineskip\vfill
\hbox to15pt{\hfill $\bigcirc$ \hfill}
}
\hfill}
\hskip-10pt plus2pt minus2pt
\hbox to30pt{\hfill
\vbox to30pt{\nointerlineskip
\hbox to25pt{\rightarrowfill} \nointerlineskip
\hbox to25pt{\leftarrowfill} \nointerlineskip
\vfill}
\hfill}
\hskip-14pt plus2pt minus2pt
\hbox to15pt{\hfill \vbox to30pt{\vskip2pt $\bigcirc$ \vfill}\hfill}
\hskip8pt plus2pt minus2pt
\hbox to30pt{\hfill
\vbox to30pt{\nointerlineskip
\hbox to25pt{\rightarrowfill} \nointerlineskip
\hbox to25pt{\leftarrowfill} \nointerlineskip
\vfill}
\hfill}
\hskip-14pt plus2pt minus2pt
\hbox to15pt{\hfill \vbox to30pt{\vskip2pt $\bigcirc$ \vfill}\hfill}
\hfil
}
\hskip80pt}


\centerline{\hfil
\hbox to 250pt{
\hfil
\hbox to15pt{\hfill \vbox to30pt{\vskip2pt${\bold E}_7\ :$ \vfill}\hfill}
\hskip10pt plus2pt minus2pt
\hbox to15pt{\hfill\vbox to30pt{\vskip2pt $\bigcirc$ \vfill}\hfill}
\hskip15pt plus2pt minus2pt
\hbox to30pt{\hfill
\vbox to30pt{
\nointerlineskip
\hbox to25pt{\rightarrowfill} \nointerlineskip
\hbox to25pt{\leftarrowfill} \nointerlineskip
\vfill}
\hfill}
\hskip-14pt plus2pt minus2pt
\hbox to15pt{\hfill \vbox to30pt{\vskip2pt $\bigcirc$ \vfill}\hfill}
\hskip8pt plus2pt minus2pt
\hbox to30pt{\hfill
\vbox to30pt{\nointerlineskip
\hbox to25pt{\rightarrowfill} \nointerlineskip
\hbox to25pt{\leftarrowfill} \nointerlineskip
\vfill}
\hfill}
\hskip-5pt plus2pt minus2pt
\hbox to 20pt{\hfill
\vbox to 30pt {\vskip2pt
\hbox to15pt{\hfill $\bigcirc$ \hfill}
\nointerlineskip\vfill
\hbox to15pt{\hfill $\uparrow\downarrow$ \hfill}
\nointerlineskip\vfill
\hbox to15pt{\hfill $\bigcirc$ \hfill}
}
\hfill}
\hskip-10pt plus2pt minus2pt
\hbox to30pt{\hfill
\vbox to30pt{\nointerlineskip
\hbox to25pt{\rightarrowfill} \nointerlineskip
\hbox to25pt{\leftarrowfill} \nointerlineskip
\vfill}
\hfill}
\hskip-14pt plus2pt minus2pt
\hbox to15pt{\hfill \vbox to30pt{\vskip2pt $\bigcirc$ \vfill}\hfill}
\hskip8pt plus2pt minus2pt
\hbox to30pt{\hfill
\vbox to30pt{\nointerlineskip
\hbox to25pt{\rightarrowfill} \nointerlineskip
\hbox to25pt{\leftarrowfill} \nointerlineskip
\vfill}
\hfill}
\hskip-14pt plus2pt minus2pt
\hbox to15pt{\hfill \vbox to30pt{\vskip2pt $\bigcirc$ \vfill}\hfill}
\hskip8pt plus2pt minus2pt
\hbox to30pt{\hfill
\vbox to30pt{\nointerlineskip
\hbox to25pt{\rightarrowfill} \nointerlineskip
\hbox to25pt{\leftarrowfill} \nointerlineskip
\vfill}
\hfill}
\hskip-14pt plus2pt minus2pt
\hbox to15pt{\hfill \vbox to30pt{\vskip2pt $\bigcirc$ \vfill}\hfill}
\hfil
}
\hskip50pt}

\centerline{\hfil
\hbox to 290pt{
\hfil
\hbox to15pt{\hfill \vbox to30pt{\vskip2pt${\bold E}_8\ :$ \vfill}\hfill}
\hskip10pt plus2pt minus2pt
\hbox to15pt{\hfill\vbox to30pt{\vskip2pt $\bigcirc$ \vfill}\hfill}
\hskip15pt plus2pt minus2pt
\hbox to30pt{\hfill
\vbox to30pt{
\nointerlineskip
\hbox to25pt{\rightarrowfill} \nointerlineskip
\hbox to25pt{\leftarrowfill} \nointerlineskip
\vfill}
\hfill}
\hskip-14pt plus2pt minus2pt
\hbox to15pt{\hfill \vbox to30pt{\vskip2pt $\bigcirc$ \vfill}\hfill}
\hskip8pt plus2pt minus2pt
\hbox to30pt{\hfill
\vbox to30pt{\nointerlineskip
\hbox to25pt{\rightarrowfill} \nointerlineskip
\hbox to25pt{\leftarrowfill} \nointerlineskip
\vfill}
\hfill}
\hskip-5pt plus2pt minus2pt
\hbox to 20pt{\hfill
\vbox to 30pt {\vskip2pt
\hbox to15pt{\hfill $\bigcirc$ \hfill}
\nointerlineskip\vfill
\hbox to15pt{\hfill $\uparrow\downarrow$ \hfill}
\nointerlineskip\vfill
\hbox to15pt{\hfill $\bigcirc$ \hfill}
}
\hfill}
\hskip-10pt plus2pt minus2pt
\hbox to30pt{\hfill
\vbox to30pt{\nointerlineskip
\hbox to25pt{\rightarrowfill} \nointerlineskip
\hbox to25pt{\leftarrowfill} \nointerlineskip
\vfill}
\hfill}
\hskip-14pt plus2pt minus2pt
\hbox to15pt{\hfill \vbox to30pt{\vskip2pt $\bigcirc$ \vfill}\hfill}
\hskip8pt plus2pt minus2pt
\hbox to30pt{\hfill
\vbox to30pt{\nointerlineskip
\hbox to25pt{\rightarrowfill} \nointerlineskip
\hbox to25pt{\leftarrowfill} \nointerlineskip
\vfill}
\hfill}
\hskip-14pt plus2pt minus2pt
\hbox to15pt{\hfill \vbox to30pt{\vskip2pt $\bigcirc$ \vfill}\hfill}
\hskip8pt plus2pt minus2pt
\hbox to30pt{\hfill
\vbox to30pt{\nointerlineskip
\hbox to25pt{\rightarrowfill} \nointerlineskip
\hbox to25pt{\leftarrowfill} \nointerlineskip
\vfill}
\hfill}
\hskip-14pt plus2pt minus2pt
\hbox to15pt{\hfill \vbox to30pt{\vskip2pt $\bigcirc$ \vfill}\hfill}
\hskip8pt plus2pt minus2pt
\hbox to30pt{\hfill
\vbox to30pt{\nointerlineskip
\hbox to25pt{\rightarrowfill} \nointerlineskip
\hbox to25pt{\leftarrowfill} \nointerlineskip
\vfill}
\hfill}
\hskip-14pt plus2pt minus2pt
\hbox to15pt{\hfill \vbox to30pt{\vskip2pt $\bigcirc$ \vfill}\hfill}
\hfil
}
\hskip10pt}


\centerline{\hfil
\hbox to 165pt{
\hfil
\hbox to15pt{\hfill \vbox to30pt{\vfill ${\bold F}_4\ :$ \vfill}\hfill}
\hskip10pt plus2pt minus2pt
\hbox to15pt{\hfill \vbox to30pt{\vfill $\bigcirc$ \vfill}\hfill}
\hskip8pt plus2pt minus2pt
\hbox to30pt{\hfill
\vbox to30pt{\vfill
\hbox to25pt{\rightarrowfill} \nointerlineskip
\hbox to25pt{\leftarrowfill}
\vfill}
\hfill}
\hskip-14pt plus2pt minus2pt
\hbox to15pt{\hfill \vbox to30pt{\vfill $\bigcirc$ \vfill}\hfill}
\hskip8pt plus2pt minus2pt
\hbox to30pt{\hfill
\vbox to30pt{\vfill \hbox to25pt{\hfill $2$ \hfill} \nointerlineskip
\hbox to25pt{\rightarrowfill} \nointerlineskip
\hbox to25pt{\leftarrowfill}
\vskip8pt}
\hfill}
\hskip-14pt plus2pt minus2pt
\hbox to15pt{\hfill \vbox to30pt{\vfill $\bigcirc$ \vfill}\hfill}
\hskip8pt plus2pt minus2pt
\hbox to30pt{\hfill
\vbox to30pt{\vfill
\hbox to25pt{\rightarrowfill} \nointerlineskip
\hbox to25pt{\leftarrowfill}
\vfill}
\hfill}
\hskip-14pt plus2pt minus2pt
\hbox to15pt{\hfill \vbox to30pt{\vfill $\bigcirc$ \vfill}\hfill}
}
\hskip135pt}


\centerline{\hfil
\hbox to 95pt{\hfil
\hbox to15pt{\hfill \vbox to30pt{\vfill ${\bold G}_2\ :$ \vfill}\hfill}
\hskip10pt plus2pt minus2pt
\hbox to15pt{\hfill \vbox to30pt{\vfill $\bigcirc$ \vfill}\hfill}
\hskip8pt plus2pt minus2pt
\hbox to30pt{\hfill
\vbox to30pt{\vfill \hbox to25pt{\hfill $3$ \hfill} \nointerlineskip
\hbox to25pt{\rightarrowfill} \nointerlineskip
\hbox to25pt{\leftarrowfill}
\vskip8pt}
\hfill}
\hskip-14pt plus2pt minus2pt
\hbox to15pt{\hfill \vbox to30pt{\vfill $\bigcirc$ \vfill}\hfill}
}
\hskip205pt}

\vskip0.1cm

\centerline{{\bf Table 1.} Calabi--Yau elliptic diagrams without
single arrows}
\centerline{(classical Dynkin diagrams).}

\newpage 

\centerline{\hbox{
$\bigcirc$
\hskip-5pt
\lower13pt\darrupdown{t_{12}}{t_{21}}
\hskip-5pt
$\bigcirc$} \hskip20pt where $t_{12}t_{21}>4$.
\hfil}

\vskip15pt


\centerline{\hbox{
$\bigcirc$
\hskip0pt
\lower13pt\darrupdown{t_{12}}{t_{21}}
\hskip-5pt
$\bigcirc$
\hskip0pt
\lower13pt\darrupdown{t_{23}}{t_{32}}
\hskip-5pt
$\bigcirc$}
\hskip20pt where $0<t_{12}t_{21}<4,\ 0<t_{23}t_{32}<4,\
t_{12}t_{21}+t_{23}t_{32}>4$ \hfil}

\vskip20pt

\centerline{\hbox{
\vbox to70pt{
\hbox to100pt{\hfill$\bigcirc$\hfill}
\vskip10pt\nointerlineskip\vfil
\hbox to100pt{\hfil
$t_{12}$ $\nearrow$\hskip-5pt$\swarrow$ $t_{21}$ \hfill
$t_{32}$ $\nwarrow$\hskip-4pt$\searrow$ $t_{23}$\hfil}
\nointerlineskip\vfil
\hbox to100pt{\hskip10pt \raise12pt\hbox{$\bigcirc$}\hfil
\darrupdown{t_{13}}{t_{31}} \hfil \raise12pt\hbox{$\bigcirc$}\hskip10pt}
}
}
\hskip20pt \raise30pt\hbox{\vbox{
\hbox{where\hfil}
\hbox{\hfil $0<t_{12}t_{21}<4,\ 0<t_{23}t_{32}<4,\
0<t_{31}t_{13}<4$,\hfil}
\hbox{$t_{12}t_{21}+t_{23}t_{32}+t_{31}t_{13}>3$}
}}
}

\vskip20pt


\centerline{
\hbox to15pt{\hfill
\vbox to 40pt{
\hbox to15pt{\hfill $\bigcirc$ \hfill}
\nointerlineskip\vfill
\hbox to15pt{\hfill $\uparrow\downarrow$ \hfill}
\nointerlineskip\vfill
\hbox to15pt{\hfill $\bigcirc$ \hfill}
\nointerlineskip}\hfill}
\hskip-5pt
\lower1pt\hbox{\hfill
\vbox to52pt {\darrup{2}
\vskip2pt
\darr}
\hskip-5pt
\hbox to15pt{\hfill
\vbox to 40pt{
\hbox to15pt{\hfill $\bigcirc$ \hfill}
\nointerlineskip\vfill
\hbox to15pt{\hfill $\uparrow\downarrow$ \hfill}
\nointerlineskip\vfill
\hbox to15pt{\hfill $\bigcirc$ \hfill}
\nointerlineskip}\hfill}
}
\hskip5pt
\hbox to15pt{\hfill
\vbox to 40pt{
\hbox to15pt{\hfill $\bigcirc$ \hfill}
\nointerlineskip\vfill
\hbox to15pt{\hfill $\uparrow\downarrow$ \hfill}
\nointerlineskip\vfill
\hbox to15pt{\hfill $\bigcirc$ \hfill}
\nointerlineskip}\hfill}
\hskip-5pt
\lower1pt\hbox{\hfill
\vbox to52pt {\darrup{2}
\vskip2pt
\darrup{2}}
\hskip-5pt
\hbox to15pt{\hfill
\vbox to 40pt{
\hbox to15pt{\hfill $\bigcirc$ \hfill}
\nointerlineskip\vfill
\hbox to15pt{\hfill $\uparrow\downarrow$ \hfill}
\nointerlineskip\vfill
\hbox to15pt{\hfill $\bigcirc$ \hfill}
\nointerlineskip}\hfill}
}
\hskip5pt
\hbox to15pt{\hfill
\vbox to 40pt{
\hbox to15pt{\hfill $\bigcirc$ \hfill}
\nointerlineskip\vfill
\hbox to15pt{\hfill $\uparrow\downarrow$ \hfill}
\nointerlineskip\vfill
\hbox to15pt{\hfill $\bigcirc$ \hfill}
\nointerlineskip}\hfill}
\hskip-5pt
\lower1pt\hbox{\hfill
\vbox to52pt {\darrup{2}
\vskip2pt
\darrdown{2}}
\hskip-5pt
\hbox to15pt{\hfill
\vbox to 40pt{
\hbox to15pt{\hfill $\bigcirc$ \hfill}
\nointerlineskip\vfill
\hbox to15pt{\hfill $\uparrow\downarrow$ \hfill}
\nointerlineskip\vfill
\hbox to15pt{\hfill $\bigcirc$ \hfill}
\nointerlineskip}\hfill}
}
\hfil}


\vskip20pt

\centerline{
\hbox to15pt{\hfill
\vbox to 40pt{
\hbox to15pt{\hfill $\bigcirc$ \hfill}
\nointerlineskip\vfill
\hbox to15pt{\hfill $\uparrow\downarrow$ \hfill}
\nointerlineskip\vfill
\hbox to15pt{\hfill $\bigcirc$ \hfill}
\nointerlineskip}\hfill}
\hskip-5pt
\lower3pt\hbox to65pt{\hfill
\vbox to50pt {
\hbox to60pt{\hfill
\vbox to12pt{
\hbox to60pt{\hfill $2$ \hfill}
\nointerlineskip
\hbox to60pt{\rightarrowfill} \nointerlineskip\vskip-2pt
\hbox to60pt{\leftarrowfill}
\vfill}\hfill}
\vfill
\hbox to60pt{\hfill
\hbox to23pt{\hfill
\vbox to12pt{\vfill
\hbox to20pt{\rightarrowfill} \nointerlineskip\vskip-2pt
\hbox to20pt{\leftarrowfill}
\vfill}\hfill}
\hskip0pt \raise4pt\hbox{$\bigcirc$} \hskip0pt
\hbox to20pt{\hfill
\vbox to12pt{\vfill
\hbox to20pt{\rightarrowfill} \nointerlineskip\vskip-2pt
\hbox to20pt{\leftarrowfill}
\vfill}\hfill}
\hfill}
}
\hfill}
\hskip-5pt
\hbox to15pt{\hfill
\vbox to 40pt{
\hbox to15pt{\hfill $\bigcirc$ \hfill}
\nointerlineskip\vfill
\hbox to15pt{\hfill $\uparrow\downarrow$ \hfill}
\nointerlineskip\vfill
\hbox to15pt{\hfill $\bigcirc$ \hfill}
\nointerlineskip}\hfill}
\hfil}

\vskip20pt

\centerline{{\bf Table 4.3.} Calabi--Yau E-sets 
without single arrows}
\centerline{(classical Lanner diagrams).}

\vskip0.5cm 

For classification of elliptic and E-subsets,  
one needs the statement:  any divisorial extremal ray 
$R$ of a 3-dimensional Calabi-Yau manifold has a curve $C\in R$ such that 
$C\cdot D(R)=-k$ where $1\le k\le 3$. This was proved by Shokurov (see  
Appendix to \cite{N13}). 

\smallpagebreak 

As a result, using Lemma 3, we get some estimate on $\rho (X)$. 
The estimate $\rho (X)\le 40$ of Theorem 5 requires additional very 
delicate considerations with elliptic and E-subsets, 
and their combinatorics.  
\enddemo

\smallpagebreak 
  
It seems, Theorem 5 is one of the strongest known results about 
the structure of Mori (or K\"ahler) cone of 3-dimensional Calabi--Yau 
manifolds. See other related results in papers of 
Miyaoka \cite{Mi}, Oguiso \cite{Og} and Wilson \cite{W1}---\cite{W3},  
and their discussion in \cite{N13}. 

It is expected that 3-dimensional Calabi--Yau manifolds 
very often have finite polyhedral Mori cone. Morrison \cite{Mor} 
conjectured that the K\"ahler cone of a Calabi--Yau manifold $X$ is 
rational finite polyhedral up to action of the automorphism group 
Aut X. In particular, the K\"ahler (and Mori) cone is 
finite polyhedral if Aut X is finite. 
It follows that a 3-dimensional Calabi--Yau manifold has 
a finite polyhedral Mori cone if its cubic intersection hypersurface 
is non-singular. We can expect that very often. 
K3 surfaces with $\rho \ge 3$ very rare have a finite polyhedral 
Mori cone (see Sect. 3). 
We expect that for Calabi--Yau 3-folds situation is opposite. 

\vskip0.5cm

Algebraic surfaces with finite polyhderal Mori cone 
are now understood well. Our knowledge about 3-folds with 
finite polyhedral Mori cone is very poor in spite of it could be very 
useful for understanding of some very important algebraic varieties, 
for example Calabi--Yau 3-folds.

\enddocument

\end

\end